\documentclass[11pt]{amsart}
\usepackage{amssymb,amsmath,amsthm}
\usepackage{epsfig}
\usepackage{graphicx}

\newtheorem{theorem}{Theorem}[section]
\newtheorem{lemma}[theorem]{Lemma}

\newtheorem{corollary}[theorem]{Corollary}

\newtheorem{example}[theorem]{Example}
\long\def\symbolfootnote[#1]#2{\begingroup%
\def\thefootnote{\fnsymbol{footnote}}\footnote[#1]{#2}\endgroup}

\newcommand\Z{{\mathbb Z}}
\begin{document}

\title{The second quandle homology of the Takasaki quandle of an odd abelian group is an exterior
 square of the group}
\date{June 1, 2010}
\author{Maciej Niebrzydowski}
\address[Maciej Niebrzydowski]{University of Louisiana at Lafayette, Department of Mathematics,
	 217 Maxim D. Doucet Hall, 1403 Johnston Street,  Lafayette, LA 70504}
\email{mniebrz@gmail.com}
\author{J\'ozef H. Przytycki}
\address[J\'ozef H. Przytycki]{The George Washington University, Department of Mathematics,
	 Monroe Hall, Room 240, 2115 G Street NW,
	 Washington, D.C. 20052; The University of Texas at Dallas, Department of Mathematics, 
	 Richardson, TX 75080}
\email{przytyck@gwu.edu}
\keywords{Takasaki quandle, quandle homology, exterior product}
\subjclass[2000]{Primary: 55N35; Secondary: 18G60, 57M25}

\thispagestyle{empty}

\begin{abstract}
We prove that if $G$ is an abelian group of odd order
then there is an isomorphism from the second quandle homology $H_2^Q(T(G))$ to
$G \wedge G$,
where $\wedge$ is the exterior product. In particular, for
$G=\Z_k^n$, $k$ odd, we have  \ $H_2^Q(T(\Z_k^n)) = \Z_k^{n(n-1)/2}$.
Nontrivial $H_2^Q(T(G))$ allows us to use 2-cocycles to construct new quandles from $T(G)$,
and to construct link invariants.
\end{abstract}

\maketitle

\section{Introduction}
Mituhisa Takasaki introduced the notion of kei (involutive quandle 
in Joyce's terminology \cite{Joy}) in 1942 \cite{Tak}.
His main example was the quandle of an abelian group, $T(G)$, with $a*b=2b-a$, which we call
Takasaki quandle. Quandle homology was first introduced in \cite{CJKLS} as a modification of 
rack homology invented by Fenn, Rourke, and Sanderson in 1995 \cite{FRS}.
We prove that if $G$ is an abelian group of odd order,
then there is an isomorphism from the second quandle homology $H_2^Q(T(G))$ to
$G \wedge G$,
where $\wedge$ is the exterior product. That is, we prove:
\begin{theorem}\label{Theorem 1.1}
$$H_2^Q(T(G))= G \wedge G,$$ 
where $G$ is an abelian group of odd order.
In particular, for
$G=\Z_k^n$, $k$ odd, we have: 
$$H_2^Q(T(\Z_k^n)) = \Z_k^{n\choose 2}.$$
\end{theorem}
This partially solves Problem 5.4 from the Ohtsuki's problem list \cite{Oht}, and is a part of our program to effectively compute homology of quandles \cite{N-P-2,N-P-3}. For other results concerning 
second homology of quandles  see \cite{Gr,E-G,Eis}.

For a quandle $(Q,*)$, we write $c\,\bar{*}\,b=a$ if $a*b=c$. We also write $a\circ c= \{b\ | \ a*b=c\}$. If
$a\circ c$ has always one element, we say that $Q$ is a quandle quasigroup. Takasaki quandle $T(G)$ is 
a quasigroup iff $G$ is a group of odd order. Then, $a\circ c= \frac{a+c}{2}$.

\section{Kernel of $\partial_2\colon C_2^Q \to C_1^Q$}
For a quandle quasigroup $Q$, the kernel of the boundary map $\partial_2\colon C_2^Q \to C_1^Q$ 
has an easy and natural description:
\begin{lemma}\label{Lemma 2.1}
Choose any element $b_0$ in a quandle quasigroup $Q$. The quotient map $C_2^Q \to 
C_2^Q/((b_0,a))$ is an isomorphism when restricted to $ker \partial_2$. Here, $((b_0,a))$ is the 
subgroup of $C_2^Q$ generated by elements of the form $(b_0,a)$ for any $a\in Q$. 
\end{lemma}
\begin{proof}
We have  $C_2^Q= C_2^{a_0} \oplus C_2^{na_0}$, where $C_2^{a_0}$ is a subgroup of $C_2^Q$ generated 
by elements $(a_0,b)$, and $C_2^{na_0}$ is the subgroup of $C_2^Q$ generated 
by all other elements of $Q^2$. Because $Q$ is a quandle quasigroup, 
$\partial_2\colon C_2^{a_0} \to C_1(Q)={\Z}Q$ is an isomorphism onto the image $\partial_2(C^Q_2)$ (which has 
 basis elements of the form $(a_0-c)$, where $c\neq a_0$). Namely, $\partial_2(a_0,b)= a_0 - a_0*b$, 
and the inverse map is given by $\partial_2^{-1}(a_0-c)= (a_0,a_0\circ c)$. The conclusion of Lemma \ref{Lemma 2.1} 
follows.
\end{proof}
\begin{corollary}
For a quasigroup quandle $Q$, we have:\\
$H^Q_2(Q)= \Z(Q^2)/((a,a),(a_0,a), im \partial_3)$.
\end{corollary}

Once again we use properties of a quandle quasigroup, to change the notation for elements of $Q^2$; namely, 
we write $[a,a*b]$ for $(a,b)$. Because $Q$ is a quandle quasigroup, $[a,c]$ is the unique 
element $(a,a\circ c)$. In this notation, the boundary map has the form: 
$$\partial_2(a,b,c)=(a,c)+ (a*c,b*c) - (a,b)- (a*b,c)=$$ 
$$[a,a*c]+ [a*c, (a*b)*c] - [a,a*b] - [a*b,(a*b)*c].$$
After substituting $x=a$, $y=(a*b)*c$, $z=a*c$, and $z'=a*b=((a*b)*c)\bar{*}(a\circ (a*c))=y\bar{*}(x\circ z)$,
we get: $[x,z]+ [z, y] - [x,z'] - [z',y]$. In the case of $Q=T(G)$, for an odd order abelian group $G$,
we have $z'=y\bar{*}(x\circ z)= z+x-y$, and consequently:
$$ \partial_2(a,b,c)= [x,z]+ [z, y] - [x,z+x-y] - [z+x-y,y].$$ 
For a quandle quasigroup, $[x,x]=0$ iff  $(x,x)=0$ and $[a_0,x]=0$ iff $(a_0,a_0\circ x)=0.$ Therefore, 
for a quandle quasigroup: 
$$H^Q_2(Q)= \Z(Q^2)/([x,x],[a_0,x],[x,z]+ [z, y] - [x,z'] - [z',y]),$$ and for 
a Takasaki quandle of an odd order abelian group $G$, we have:
$$H^Q_2(T(G))= 
\Z(G\times G)/([x,x],[0,x],[x,z]+ [z, y] - [x,z+x-y] - [z+x-y,y]).$$ 
We show in Main Lemma that this group is in fact $G\wedge G$.

\section{Main Lemma}
\begin{lemma} 
Let $G$ be an abelian group of odd order. Then:\\
$\Z(G\times G)/([x,z]+ [z,y] - [x,z+x-y] - [z+x-y,y], [x,x], [0,x]) \ = \ G\wedge G$. 
\end{lemma}

\begin{proof} 
If we add the bilinearity relation $[a,b+c]-[a,b]- [a,c]$ to relations on the left side of the equation,
then the long relation automatically holds, and we get an epimorphism $H^Q_2(T(G)) \to G\wedge G$.
To prove isomorphism, we show that the bilinearity follows from our relations:\\
Step 1: (antisymmetry) $[x,z]=-[z,x]$.\\
To demonstrate this relation, assume that $y=x\circ z = \frac{x+z}{2}$ (thus, $z+x-y= y$). Then we get:\\
$[x,z]+[z,y] = [x,y]+[y,y] = [x,y]$. Exchanging the roles of $x$ and $z$ gives:\\
$[z,x]+[x,y] = [z,y]+[y,y] = [z,y]$.\\
Adding these equalities sidewise we get:\\
$[x,z] + [z,x] = 2[y,y] = 0$, as needed. In particular, $[x,0]=0$.\\
Step 2: $[x,z]=[x, z+x]$.\\
We obtain this identity by substituting $y=0$ in $[x,z]+ [z,y]=[x,z+x-y]+ [z+x-y,y]$, and applying the last relation from Step 1. 
In particular, it follows that:
$[x,-z]= [x+z,-z]= [x+z,x]=[z,x]=-[x,z]$.\\
Step 3. We are now ready for the last calculation showing the relation\\
$[w,v+z]= [w,v] + [w,z]$:

\subsection{Final calculation for bilinearity}
\ \\
In the identity $[x,z]+[z,y] = [x,z-y] + [x+z,y]$ obtained from Step 2, we substitute 
$v=y-z$ (so $y=z+v$) to get:
$$[x,z]+ [z,z+v] = [x,-v] + [x+z, z+v],$$
and by already proven identities: $[z,z+v] =[z,v]$, $[x,v]= -[v,x]$, $[x,-v]=-[x,v]$, we obtain:
$$[z+x, z+v]= [x,z] + [z,v] + [x,v].$$ 
Then, by exchanging the roles of $x$ and $z$:
$$[x+z, x+v]= [z,x] + [x,v] + [z,v].$$ Taking the difference gives:
$$[z+x, z+v] - [z+x, x+v] = 2[x,z].$$ Let $w=x+z$ (i.e., $x=w-z$). Then:
$$[w,z+v] - [w, w-z+v]= 2[w-z,z],$$ or equivalently:
$$[w,z+v] + [w,z-v]= 2 [w,z].$$
Exchanging the roles of $v$ and $z$ gives: 
$$[w,v+z] + [w,v-z]= 2[w,v],$$ and by adding the last two equalities sidewise, we obtain:
$$2[w,v+z] = 2[w,z]+ 2[w,v].$$
Because $G$ has an odd order, 2 is not a zero divisor, and we get:
$$[w,v+z] - [w,v] - [w,z]=0, $$ as needed.
Finally, we use the fact (constructive definition) that for $G$ of odd order:
$$G\wedge G = (G\otimes G)/([x,y]+[y,x])= 
\Z(G\times G)/([x,y+z]-[x,y] - [x,z], [x,y]+[y,x]).$$ 
This completes the proof of Main Lemma and Theorem \ref{Theorem 1.1}.
\end{proof}

\section{Conclusions}
Nontrivial $H_2^Q(T(G))$ allows us to use 2-cocycles to construct new quandles from $T(G)$ 
(in fact involutive quandles, as observed by M. Saito). 
We illustrate it with an example of $G=\Z_p\oplus \Z_p$, in which case $H_2^Q(T(\Z_p^2))=\Z_p$.

\begin{example} By the universal coefficient theorem
for cohomology,\\ 
$H^2_Q(T(\Z_p^2))= \Z_p$, and the generating cocycle 
can be written as $$\phi=(e^*_1\wedge e^*_2)\mathcal P,$$ where $\mathcal P\colon C_2^Q(\Z_p^2) \to \Z_p^2 \wedge \Z_p^2$ is 
the quotient map from the Main Lemma, 
$e_1=(1,0),e_2=(0,1)$ is a basis of $\Z_p^2$, 
and $e^*_1$, $e^*_2$ is the dual basis (i.e., basis of $Hom(\Z_p^2 \to \Z_p)$). Furthermore, $e^*_1\wedge e^*_2$ is 
the induced map from $\Z_p^2 \wedge \Z_p^2$ to $\Z_p$. In particular, $\phi(e_1,e_2)=1=-\phi(e_2,e_1)$. By \cite{CKS}, 
one can construct a new quandle $E(T(\Z_p^2),\Z_p,\phi)$ on the set $\Z_p^2\times \Z_p$, with quandle operation 
$(x_1,a_1) * (x_2,a_2)= (x_1*x_2, a_1 + \phi(x_1,x_2))$. It is called the central extension of $T(\Z_p^2)$ by $\Z_p$ using the cocycle $\phi$.  
This kei is not isomorphic to $T(\Z_p^3)$, as it is not a quandle quasigroup. We checked by GAP \cite{GAP} that 
$H^Q_2(E(T(\Z_3^2),\Z_3,\phi))=0$.
\end{example}
\begin{example} 
There are five different groups of order 27. We consider their core quandles 
(recall that $core(G)$ is an involutive quandle with $g*h=hg^{-1}h$; the same as Takasaki quandle in abelian case).
Their second quandle homology is given below:\\
(1) $G_1=\Z_{27}$; $H^Q_2(T(G_1))= 0$.\\
(2) $G_2=\Z_{3}\oplus \Z_9$; $H^Q_2(T(G_2))= \Z_3$.\\
(3) $G_3=\Z_{3}^3$; $H^Q_2(T(G_3))= \Z_3^3$.\\
(4) $G_4=\{s,t\ | \  s^9=t^3=1,\ st=ts^4\}$; $H^Q_2(core(G_4))=\Z_3$.\\
(5) $G_5=\{x,y,z\ | \  x^3 = y^3 = z^3 = 1, yz = zyx, xy = yx, xz = zx \}$; $H^Q_2(core(G_5))= \Z_3^3$.\\
$G_5$ is the modulo 3 Heisenberg group. It is also a quotient  of  the Burnside group of exponent 3 ( $B(3,3)=  \{x,y,z\ | \  w^3 =1,\ for \ any\ word \ w\}$,
and $core(G_5)$ is a commutative kei \cite{N-P-1}.

We challenge the reader to use quandle homology to distinguish the quandle $T(G_2)$ from $core(G_4)$, 
and  $T(G_3)$ from $core(G_5)$.

\end{example}

\begin{figure}
\begin{center}
\includegraphics[height=6cm]{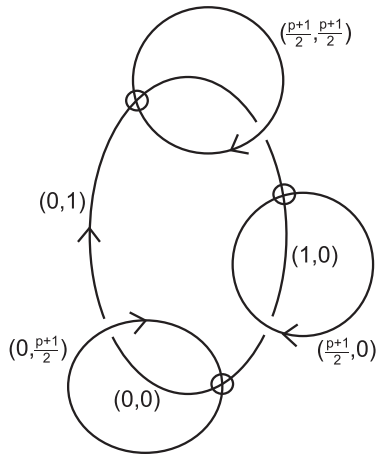}
\caption{2-cycle realized by a virtual link diagram\label{przyklad} }
\end{center}
\end{figure}

2-cycles in quandle homology correspond to colored virtual knot diagrams.
Each positive crossing represents a pair $(x,y)\in C^Q_2(X)$, where $x\in X$ is the color of an under-arc away 
from which points the normal of the over-arc labeled by $y\in X$. 
In the case of negative crossing, we write $-(x,y)$.
The sum of such 2-chains taken over all crossings of the diagram forms a 2-cycle (see \cite{Gr,CKS} for more details).
An example of a virtual link diagram realizing our generator $e_1\wedge e_2 \in H^Q_2(T(\Z_p^2))$ is shown in the Figure \ref{przyklad}.
It follows that this virtual link is nontrivial.

In the sequel paper we extend our results to all finitely generated Takasaki quandles and some 
other Alexander quandles. In particular, we complete Greene's calculations \cite{Gr} by showing that 
$tor H^Q_2(T(\Z_{4k}))= \Z_2^2$.

\section{Acknowledgements}
M.~Niebrzydowski was partially supported by the Louisiana Board of Regents grant (\# LEQSF(2008-11)-RD-A-30). 
J.~H.~Przytycki was partially supported by the NSA grant (\# H98230-08-1-0033), 
by the Polish Scientific Grant: Nr. N N201387034, by the GWU REF grant,
 and by the CCAS/UFF award.

\end{document}